\documentclass[12pt]{article}
\usepackage{amscd}
\usepackage{mathrsfs}
\usepackage{amsfonts}
\usepackage{graphicx}
\usepackage{bookmark}
\usepackage{amsmath,amscd,stmaryrd,amssymb,array, amsfonts,amsmath,amssymb}
\usepackage{enumitem}
\usepackage{color}
\usepackage{CJK}
\usepackage{appendix}
\setenumerate[1]{itemsep=0pt,partopsep=0pt,parsep=\parskip,topsep=5pt}
\setitemize[1]{itemsep=0pt,partopsep=0pt,parsep=\parskip,topsep=0pt}
\setdescription{itemsep=0pt,partopsep=0pt,parsep=\parskip,topsep=5pt}

\numberwithin{figure}{section}
 \numberwithin{equation}{section}

\newtheorem{theorem}{Theorem}[section]
\newtheorem{proposition}[theorem]{Proposition}
\newtheorem{definition}[theorem]{Definition}

\newtheorem{lemma}[theorem]{Lemma}
\newtheorem{remark}[theorem]{Remark}

\newcommand{\cC}{{\mathcal C}}

\newcommand{\cL}{{\mathcal L}}

\newcommand{\sL}{{\mathscr L}}

\def\be{\begin{equation}}
\def\ee{\end{equation}}
\def\bes{\begin{equation*}}
\def\ees{\end{equation*}}
\def\bsp{\begin{split}}
\def\esp{\end{split}}

\def\ba{\begin{array}}
\def\ea{\end{array}}
\def\benu{\begin{enumerate}}
\def\eenu{\end{enumerate}}
\def\bt{\begin{theorem}}
\def\et{\end{theorem}}
\def\bp{\begin{proposition}}
\def\ep{\end{proposition}}
\def\bl{\begin{lemma}}
\def\el{\end{lemma}}
\def\br{\begin{remark}}
\def\er{\end{remark}}
\def\bd{\begin{definition}}
\def\ed{\end{definition}}

\def\b{\beta}
\def\De{\Delta}
\def\de{\delta}
\def\pa{\partial}
\def\nab{\nabla}
\def\lam{\lambda}
\def\Lam{\Lambda}

\def\ve{\varepsilon}
\def\vth{\vartheta}
\def\sig{\sigma}

\def\gam{\gamma}

\def\k{\kappa}

\def\a{\alpha}
\def\W{\Omega}

\def\.{\cdot}

\def\R{\mathbb{R}}

\def\A{\forall}

\def\Cap{\bigcap}
\def\ln{\mbox{ln\,}}
\def\ra{\rightarrow}

\def\~{\tilde}
\def\8{\infty}
\def\X{\times}

\def\mb{\mbox}

\def\Hs{\hspace{1cm}}\def\hs{\hspace{0.5cm}}
\def\Vs{\vskip8pt}\def\vs{\vskip4pt}

\def\({\left(}\def\){\right)}

\begin{document}

\begin{center}
{\bf\Large {Global Estimates  and Regularity of Retarded\\ Parabolic Equations with  Fast-growing \\[0.2cm] Nonlinearities \footnote{This work was supported by the National Natural Science Foundation of China [11871368]. }}}
\end{center}

\vs\centerline{Desheng  Li}
\begin{center}
{\footnotesize
{School of Mathematics,  Tianjin University\\
          Tianjin 300072,  China\\

{\em E-mail}:  lidsmath@tju.edu.cn (D. Li), }}
\end{center}

{\footnotesize
\noindent{\bf Abstract.} This paper is  concerned with  global estimates  and regularity of solutions for  the initial value problem of the retarded parabolic equation
$$\frac{\pa u}{\pa t}-\De u=f(x,u)+g(u(x,t-r_1(t)),\cdots,u(x,t-r_m(t)))+h(x,t)$$ in a bounded domain $\W\subset \R^n$ with fast-growing  nonlinearities and a dissipative structure, which is associated with the homogeneous Dirichlet boundary condition. Our results reveal some deeper inherent connections between dissipative structures and the regularity of solutions for such problems.
 \Vs
\Vs \noindent{\bf 2010 MSC:}  35B40, 35B41, 35B65, 35K20, 35K58
}

\section{Introduction}

This paper is basically concerned with global estimates and regularity of solutions for the initial value problem of the following retarded parabolic equation:
\be\label{e1.1}\begin{split}
\frac{\pa u}{\pa t}-\Delta u=&f(u)+g(u(x,t-r_1(t)),\cdots,u(x,t-r_m(t)))
+h(x,t)
\end{split}\ee
with fast growing nonlinearities  in a bounded domain $\W\subset \R^n$, which is  associated with the homogeneous Dirichlet boundary condition:
\be\label{e1.2}
u(x,t)=0,\Hs x\in\pa\W,
\ee
 where  $f\in C^1(\R), \,g\in C^1(\R^m)$, and   $h$  is a measurable  function on $\W\X \R$. The delay functions $r_i$ belongs to $C(\R;[0,r])$ ($1\leq i\leq m$) for some $0\leq r<\8$.
This type of problems have  been  studied by many authors in the past decades; see e.g. \cite{Chue,Hale1,Hale2, KL, MR,SC, TW1,TW2,WK, WW,Wu, ZS}. However, due to technical difficulties induced by time lags in the equations, we find   that most of the existing works mainly  focus  on the case where the delay terms  have at most  sublinear nonlinearities.  

In this present work we are interested in  the case where both the dominant term  $f$ and the delay term $g$ in \eqref{e1.1} may have arbitrary polynomial growth rates. Since we make no restrictions on the dimension of the space,
the investigation of qualitative properties  such as the  global existence and  uniqueness, the regularity, and particularly the dynamics of the equation in functional spaces with higher regularities is in fact never an easy task even if we come back to the situation of the non-delayed case; see e.g. \cite{Hale1,Robin,Sell,Tem, Vish1,Vish2} etc.
Our main  purpose is to carry out a systematic study on the aforementioned problem under a typical  dissipative type condition on  $f$.

Specifically, let $f$ and $g$ satisfy the following structure conditions:
\benu
\item[{(F0)}]  There exist $\Lam,N>0$ and   $\gam>1$ such that
$$
f(s)s\leq -\Lam |s|^{{\gam}+1}+N,\Hs \A\,s\in\R.
$$
\item[{(F1)}]  There exists  ${\a}\geq 1$ and $a_0>0$ such that
$$
|f'(s)|\leq a_0 (|s|^{{\a}-1}+1),\Hs \A\,s\in\R.
$$
\item[(G1)] There exist $1\leq\b<\gam$ and $b_0>0$ such that
$$
 |\nab g(v)|\leq b_0 (|v|^{{\b}-1}+1),\Hs \A\,v=(v_1,\cdots,v_m)\in\R^m.
 $$
 \eenu
For simplicity in writing, we assign the initial time $\tau=0$ for \eqref{e1.1}-\eqref{e1.2} and write $u(t)=u(\.,t)$ for a solution $u(x,t)$   (in the distribution sense) of the problem with initial value
\be\label{eic}
u(\.,\,s)=\phi(s),\Hs s\in[-r,0]
\ee
Set
\be\label{e1.p*}
p_*={\b}({\gam}-1)/({\gam}-{\b}),\hs q_*=\max(p_*,2\a,2\b).\ee
Denote $h(t)=h(\.,t)$. 

\vs Given  $q_*<q\leq \8$, suppose $h\in L^\8(\R;L^{q_\gam/\gam}(\W))$, where $q_\gam=q-1+{\gam}$. We show that \eqref{e1.1}-\eqref{eic} has a unique global solution
$$ u\in C\([-r,\8);V_1\)\Cap L^\8\((-r,\8);V_1\)\Cap L^\8\((-r,\8);L^q(\W)\)$$  for any initial value function 
$$\phi\in C\([-r,0];{V_1}\)\Cap L^\8\((-r,0);L^q(\W)\),$$ where $V_1=H^1_0(\W)$.
More importantly, we establish exponential decay estimates in $V_1\Cap L^q(\W)$ in case $q<\8$ and uniform  boundedness  estimates in $V_1\Cap L^q(\W)$ in case $q=\8$.

Concerning the regularity properties  of $u$ and the  dissipativity of the problem in functional spaces with higher regularities, we have  the following interesting results:
\vs
Assume $r_i\in C^1(\R)$ $(1\leq i\leq m)$, and suppose that
\be\label{e1.9}\left\{\ba{ll}
 h\in L^\8(\R\X\W)\Cap L^\8(\R;H^1(\W)),\\[1ex]
  h'\in L^2((a,b);H),\hs \A\,-\8<a<b<\8,\ea\right.\ee
 Then 
 \be\label{e1.10}\left\{\ba{ll}
 u\in C([-r,\8);V_2)\Cap L^\8((-r,\8);V_2),\\[1ex]
  u'\in L^2((0,T);V_1)\Cap C([0,T];H),\hs \A\,0<T<\8,\ea\right.\ee
 provided that
 $$\phi\in  C\([-r,0];{V_2}\)\Cap L^\8\((-r,0)\X\W\)$$ and  $\phi'\in L^2((-r,0);H)$, where
 $V_2=H^2(\W)\Cap H^1_0(\W)$.
   \vs
   A particular but important case in applications is the one  where the time lags  $r_i$ ($1\leq i\leq m$) are separated, i.e. the function $g$ takes the form
\be\label{e1.8}
g(v)=g_1(v_1)+g_2(v_2)+\cdots+g_m(v_m),\Hs v\in\R^m.
\ee
In such a case we can show that if 
$$\phi\in C\([-r,0];{V_2}\)\Cap L^\8\((-r,0);L^q(\W)\),\hs \int_{-r}^0\int_\W|\phi|^{q-2}|\nab \phi|^2dx\,dt<\8$$ for some $q_*<q<\8$ and
$
\phi'\in L^2((-r,0);H)$, 
 then the regularity results in \eqref{e1.10} remain true.  What is more,  for each fixed $q>q_*$, there exist  positive constants $\~M_2,\lam_2$ and $\~\rho_2$ such that 
\bes\label{e:3.46}
 \begin{split}
|\De u|^2\leq \~M_2\(\|\phi\|_{C\([-r,0];{V_2}\)}^2+\|\phi\|_{L^\8\((-r,0);L^q(\W)\)}^q+\int_\W|\phi|^{q-2}|\nab \phi|^2dx\)e^{-\lam_2t}+\~\rho_2\,.
\end{split}\ees
\vs

It is interesting to note that for  $1\leq q\leq q_*$ (and particularly for $1\leq q\leq p_*$), the global existence of \eqref{e1.1}-\eqref{eic} remains an open question. To the best of our knowledge, this is the case even if for the non-delayed parabolic equations with nonlinearities as considered here.

This paper is organized as follows. In Section 2 we do some preliminaries, and in   Section 3 we establish global $L^q$ and $H^1$ estimates. Section 4 is devoted to global $H^2$ estimates. In  Section 5 we state our main results on the global  existence and uniqueness, and the regularity of solutions of the problem.

\section{Preliminaries}
This section is concerned with some preliminaries. We first recall some fundamental inequalities. Then we give an abstract form of the initial value problem of system \eqref{e1.1}-\eqref{e1.2}.
\subsection{Some fundamental  inequalities}
First, we have the following easy facts.
\bl\label{l:2.0a}Let $p>0$. Then for any $a_1,a_2,\cdots,a_m\geq 0$,
$$
(a_1+a_2+\cdots+a_m)^p\leq m^p(a_1^p+a_2^p+\cdots+a_m^p).
$$
\el

\bl\label{l:2.0}Let $p,q>1$, and $1/p+1/q=1$. Then for any $a,b,\ve>0$,
$$
ab\leq \ve a^p+\ve^{-q/p}b^q.
$$
\el
{\bf Proof.} This is a simple consequence of the classical  Young's inequality.  Indeed,
$$
ab\leq \ve a^p+\frac{1}{qp^{q/p}}\ve^{-q/p}\,b^q\leq \ve a^p+\ve^{-q/p}b^q.
$$
This  is precisely what we desired.  $\Box$

\Vs
Using the H\"{o}lder's inequality and Lemma \ref{l:2.0}, one trivially verifies the validity of the lemma below.

\bl\label{l:2.0c}Let $0<p<q<\8$, and let $\W\subset \R^n$ be a bounded domain. Then for any  $u\in L^q(\W)$, $v\in L^{(q-p)/p}(\W)$ and $\ve>0$,
$$
\int_\W|u|^p|v|dx\leq \ve |u|_q^q+\ve^{-p/(q-p)}\left|v\right|_{(q-p)/p}^{(q-p)/p}.
$$
In particular,
$$
\int_\W|u|^pdx\leq \ve |u|_q^q+\ve^{-p/(q-p)}|\W|.
$$
Here and below $|\W|$ denotes the Lebesgue measure of $\W$.
\el

The following two    retarded integral inequalities given in a recent paper \cite{LLJ} by Li {\em et al.} will play a crucial role in our argument.

Let $E$ be a bounded nonnegative measurable function on $Q:=(\R^+)^2$ satisfying
\be\label{e1.13}\ba{ll}
\lim_{t\ra +\8}E(t+s,s)=0\mb{ uniformly w.r.t.  $s\in\R^+$}, \ea
\ee and let $K$ be a  nonnegative measurable function on $Q$ with
\be\label{e1.14}\ba{ll}
\kappa:=\sup_{t\geq 0}\int_0^t K(t,s)ds<\8.\ea
\ee
Given $r\geq 0$, denote  ${\cC}$ the space $C([-r,0])$ equipped with the usual sup-norm $\|\.\|_\cC$.
Consider the retarded integral inequality
\be\label{e1.3}\ba{ll}
y(t)\leq &E(t,\tau)\|y_\tau\|_\cC+\int_\tau^t K(t,s)\|y_s\|_\cC ds+\rho,\Hs\A\,t\geq\tau\geq 0,
\ea
\ee
where $\rho\geq 0$ is a constant, $y$ is a nonnegative continuous function, and $y_t$ denotes the {\em element} in ${\cC}$,
\be\label{e1.15}
y_t(s)=y(t+s),\Hs s\in [-r,0].
\ee

Denote $\sL_r(E;K;\rho)$ the solution set of \eqref{e1.3}, namely,
$$\sL_r(E;K;\rho)=\{y\in C([-r,\8)):\,\,\,y\geq0\mb{ and satisfies } \eqref{e1.3}\}. $$

\bl\label{l:2.4}\cite[Theorem 1.3]{LLJ} 
\,  The the following two assertions hold.
  \benu
  \item[$(1)$] If $\k<1$ then for any  $R,\ve>0$, there exists $T>0$ such that
  \be\label{e:t2.2}
  y(t)<\mu \rho+\ve,\Hs t>T
  \ee
  for all   $y\in \sL_r(E;K;\rho)$ with $\|y_0\|_\cC\leq R$, where
\be\label{emu}
 \mu =1/(1-\k).\ee
\vs\item[$(2)$] If  $\k<1/(1+\vartheta)$, where $\vth=\sup_{t\geq s\geq0}E(t,s)$, then there exist   $M,\lam>0$ $($independent of $\rho$$)$ such that
\be\label{e:gi}
{y(t)}\leq M\|y_0\|_\cC e^{-\lam t}+\eta\rho,\Hs t\geq0
\ee
for all  $y\in \sL_r(E;K;\rho)$, where
\be\label{ec}
\eta=({\mu+1})/{(1-\k c)},\hs c =\max\(\vartheta /(1-\kappa),\,1\).
\ee
\eenu
\el

\br\label{r1.1}If $\k<1/(1+\vartheta)$ then one  trivially verifies that
$
\k c <1.
$

\er

\br\label{r:2.5}
In most examples from applications, the function $E(t,s)$ in \eqref{e1.3} is an exponential function,
$$
E(t,s)=M_0e^{-\lam_0(t-s)},
$$
where $M_0$ and $\lam_0$ are positive constants. In such a case we infer from \cite[Remark 2.2]{LLJ} that the constants $M$ and $\lam$ in \ref{e:gi} can be taken as
$$
M=c\sqrt{2/(1+\k c)},\hs \lam=\frac{\ln2-\ln\({1+\k c}\)}{2\(M_1+r\lam_0\)}\,\lam_0,
$$
{where }$M_1=\max\(\ln(M_0\eta),\,\ln\(\frac{2M_0}{1-\k c}\)\).$
In particular, if $r=0$ then we have
$$\lam=\xi\lam_0,\hs \xi=\frac{\ln2-\ln\({1+\k c}\)}{2M_1}.$$
\er
%

 \bl\label{l:2.3}\cite[Lemma 2.1]{LLJ} Suppose $\k<1$.
Let $y$ be a nonnegative continuous function on $[-r,T)$ $(0<T\leq\8)$ satisfying
  \be\label{e:3.1e}
y(t)\leq E(t,0)\|y_0\|_\cC+\int_0^t K(t,s)\|y_s\|_\cC\, ds+\rho,\Hs0\leq t<T.
\ee
  Then
\be\label{e:3.4'}
y({t})\leq (c +1) (\|y_0\|_\cC+1)+\mu\rho,\Hs t\in[0,T),
\ee
where $\mu$ and  $c$ are the constants defined in Theorem \ref{t:3.1}.
\el

\subsection{Abstract form of problem \eqref{e1.1}-\eqref{e1.2}}
Let $H=L^2(\Omega)$, $V_1=H^1_0(\Omega)$, and $V_2=H^2(\Omega)\Cap H^1_0(\Omega)$. Denote by $(\.,\.)$ and  $|\.|$ the inner product and  norm on $L^2(\Omega)$, respectively, and define the norms $||\.||_1$ on  $V_1$ and $||\.||_2$ on  $V_2$ as follows:
$$||u||_1=|\nabla u| \,\,\,(u\in V_1),\hs ||u||_2=|\De u| \,\,\,(u\in V_2).$$
By the basic theory on fractional powers of spaces (see e.g. \cite[Chap. 1.4]{Henry}), the  $\|\.\|_i$ is equivalent to the usual norm on $V_i$.
\vs
We will also use $|\.|_q$ to denote the norm of $L^q(\Omega)\ (0<q\leq\8)$.
\vs
Denote   $A$  the operator $-\Delta$ subjects to the homogeneous Dirichlet boundary condition, and    $\mu_k$ $({k\geq 1})$ the distinct  eigenvalues of $A$,
$$
0<\mu_1<\mu_2<\cdots<\mu_k<\cdots.
$$

Let $X$ be a Banach space $X$. Given $r> 0$, denote $\cC_X$ and $\cL^\8_X$ the spaces $C([-r,0],X)$ and $L^\8([-r,0],X)$ equipped with the usual norms:
$$\|\.\|_{\cC_X}=\|\.\|_{C([-r,0],X)},\hs \|\.\|_{\cL^\8_X}=\|\.\|_{L^\8([-r,0],X)},$$ respectively. 
To deal with differential equations with and without delays in a uniform manner,  we also assign
$$
\cC_X=\cL^\8_X=X,\hs\mb{if }\,r=0.
$$
In case $X=L^q(\W)$ ($0<q\leq \8$), we will  simply write
$$\cC_{L^q(\W)}=\cC_q,\hs \cL^\8_{L^q(\W)}=\cL^\8_q.$$

\vs  The {\em lift}  of  a function 
$u\in L^\8([\tau-r,T),X)$ ($T>\tau$) in $\cL^\8_X$ is defined to be a mapping $u_t$ from $[\tau,T)$ to $\cL^\8_X$,
$$
u_t(s)=u(t+s),\Hs s\in[-r,0].
$$

Now we define informally a mapping $G(t,\phi)$ in $\R\X\cC_{V_1}$ as below:
$$
G(t,\phi)=g(\phi(-r_1(t)),\cdots,\phi(-r_m(t))),\Hs (t,\phi)\in D(G)\subset \R\X\cC_{V_1}.
$$
Set $h(t)=h(\.,t)$. Then problem \eqref{e1.1}-\eqref{e1.2} can be put into an abstract form:
 \be\label{eae}\frac{du}{dt}+Au=f(u)+G(t,u_t)+h(t).\ee
 Since equation \eqref{eae} is nonautonomous, one has to take into account the initial time when considering its initial value problem. Hence the initial value problem of the equation  generally reads as
  \be\label{e:ivp1}\left\{\ba{ll}\frac{du}{dt}+Au=f(u)+G(t,u_t)+h(t),\hs t\geq\tau,\\[1ex]
  u(\tau+s)=\phi(s),\hs s\in[-r,0],\ea\right.\ee
 where $\phi\in\cC_{V_1}$, and $\tau\in\R$ denotes the initial time. Rewriting  $t-\tau$ as $t$, one obtains an equivalent form of \eqref{e:ivp1}:
 \be\label{e:ivp2}\left\{\ba{ll}\frac{du}{dt}+Au=f(u)+G^\tau(t,u_t)+h^\tau(t),\hs t\geq 0,\\[1ex]
 u(s)=\phi(s),\hs s\in[-r,0],\ea\right.\ee
where
$$
G^\tau(t,\.)=G(t+\tau,\.),\hs h^\tau(t)=h(t+\tau).
$$

Given $\tau\in\R$ and $\phi\in\cC_{V_1}$, denote  $u(t;\tau,\phi)$ the solution $u$  of \eqref{e:ivp2} (if exists) in the distribution sense on a maximal existence interval $[-r,T_\phi)$ ($T_\phi>0$). For convenience, we call the lift $u_t$ of  $u$ the {\em solution curve} of \eqref{e:ivp2} in $\cC_{V_1}$ with initial value $u_\tau=\phi$, denoted hereafter by $u_t(\tau,\phi)$.

\section{$L^q$ and $H^1$ Decay Estimates}
In this section we establish some decay estimates for solutions of the initial value problem \eqref{e:ivp2}, which in turn  imply the regularities of the solutions in appropriate functional spaces.

For simplicity, we only consider the case where the initial time $\tau$ in \eqref{e:ivp2} equals $0$. One easily sees that  all the estimates given below for $u(t)=u(t;0,\phi)$  hold true for solutions $u(t;\tau,\phi)$ of \eqref{e:ivp2} in a uniform manner with respect to $\tau\in\R$.

\vs
It should be pointed out that many calculations leading to the estimates are not reasonable because a solution $u$ of \eqref{e:ivp2} in the distribution sense  may not be sufficiently regular. For instance, in general it remains unknown whether the $L^q$-norm $|u(t)|_q$  is a continuous function in $t$ for sufficiently large $q$. Hence in the proof of Theorem \ref{t:3.1} below,  Lemmas \ref{t:3.1} and \ref{l:2.3} can not be directly applied
to  $y(t):=|u(t)|_q$ to derive decay estimates.
However, they can be justified by considering appropriate  approximations $u_k$ of $u$ as follows.

Let $w_j(j=1,2,\cdots)$ be an othorgonal basis of $H=L^2(\W)$ consisting of  eigenvectors of $A$. Given  $\phi\in \cC_{V_1}\Cap \cL^\8_q$, pick a sequence $\phi_k=\sum_{j=1}^ka_j(t)w_j$ and a sequence $h_k=\sum_{j=1}^kb_j(t)w_j$ with sufficiently smooth coefficients $a_j(t)$ and $b_j(t)$ such that $\phi_k\ra \phi$ and $h_k\ra h$ in appropriate topologies of $\cC_{V_1}\Cap \cL^\8_q$ and $L^\8(\R,L^q(\W))$, respectively.  For each $k$, let
$$
u_k(t)=\sum_{j=1}^kc_{kj}(t)w_j
$$
be a Galerkin approximation of \eqref{e:ivp2} which  solves the following  system:
\be\label{e:ivpa}\left\{\ba{ll}\(\frac{du_k}{dt}+Au_k,w_j\)=\(f(u_k)+G(t,u_{k,t})+h_k(t),w_j\),\hs 1\leq j\leq k,\\[1ex]
 u_k(s)=\phi_k(s),\hs s\in[-r,0].\ea\right.\ee
 Here and below $(\.,\.)$ denotes the inner product of $H=L^2(\W)$.
By the basic theory on  ODEs we know that $u_k$ is sufficiently regular, so that all the calculations
can be performed rigorously on $u_k$. As a result, the estimates in the theorems below   remain valid for $u_k$. Passing to the limit one immediately concludes that these estimates hold true for $u$.

Hence in what follows we always suppose  that both the initial value function $\phi$ in \eqref{e:ivp2} and the  solution $u$ of \eqref{e:ivp2} are  sufficiently regular when performing mathematical  calculations.
\subsection{Decay Estimates in $L^q(\W)$}
We first observe that by (F1) and  (G1) one has
\be\label{e1.4}
|f(s)|\leq a_1 (|s|^{{\a}}+1),\Hs \A\,s\in\R,
\ee
and
\be\label{e1.5}
|g(v)|\leq b_1(|v|^{{\b}}+1),\Hs \A\,v\in\R^m.
\ee

Since $ \|v\|_*=\max_{1\leq i\leq m}|v_i|$ is a norm in $\R^m$ and  all the norms in $\R^m$ are equivalent, by (G1) and \eqref{e1.5} we also have
\be\label{e1.6}
|\nab g(v)|\leq b_0'(\|v\|_*^{{\b}-1}+1),\Hs \A\,v\in\R^m,
\ee
and
\be\label{e1.7}
|g(v)|\leq b_1'(\|v\|_*^{{\b}}+1),\Hs \A\,v\in\R^m.
\ee

Let $p_*$ and $q_*$ be given as in \eqref{e1.p*},  and write
$$q_\gam=q-1+{\gam},\Hs \A\,q\geq 1.
$$
Denote  $u(t;\phi)$  the solution of \eqref{e:ivp2} with $\tau=0$ and initial value $\phi$. Our
first result is summarized in the following theorem.

 \bt\label{t:3.1} Let $q_*< q<\8$. Suppose
$h\in L^\8(\R;L^{q_\gam/\gam}(\W))$.  Then  for each  $\phi\in \cC_{V_1}\Cap \cL^\8_q$, $u(t;\phi)$ is globally defined for $t\geq 0$. Furthermore, there exist $M,\lam_q,\rho_q>0$, where $M$ is independent of $q$, such that
\be\label{e:3.1}
|u(t;\phi)|_{q}^{q}\leq Me^{-\lam_q t}\|\phi\|_{\cL^\8_q}^{q}+\rho_q,\Hs t\geq 0.
\ee
   \et
{\bf Proof.}  Let  $[-r,T_\phi)$ be the  maximal existence interval of $u=u(t;\phi)$.
Taking the inner product of both sides of \eqref{eae} with $|u|^{q-2}u$ for $q\geq 2$, we obtain by (F0) and \eqref{e1.5} that
\be\label{e4.2}
\begin{split}
&\hs\,\,\frac{1}{q}\frac{d}{dt}|u|_{q}^{q}+(q-1)\int_\W|u|^{q-2}|\nab u|^2dx\\[1ex]
&\leq\int_\W|u|^{q-2}u f(u)dx+\sum_{i=1}^m\int_\W|u|^{q-1}|G(t,u_t)|dx+\int_\W |u|^{q-1}|h(t)| dx\\
&\leq -\Lam|u|_{{q_\gam}}^{{q_\gam}}+N|u|_{q-2}^{q-2}+b_1' \int_\W|u|^{q-1}(\,\max_{1\leq i\leq m}|u(t-r_i)|^{{\b}}+1)dx\\&\hs+ \int_\W |u|^{q-1}|h(t)| dx.
\end{split}
\ee
Using the H\"{o}lder's inequality and  Lemma \ref{l:2.0} we deduce  that
\begin{equation*}
\begin{split}
&\hs \int_\W|u|^{q-1}(\,\max_{1\leq i\leq m}|u(t-r_i)|^{{\b}})dx\\
&\leq \int_\W|u|^{q-1}(\,\sum_{1\leq i\leq m}|u(t-r_i)|^{{\b}})dx\leq \sum_{1\leq i\leq m}|u|^{q-1}_\sig\,|u(t-r_i)|_{q}^{{\b}}\\
&\leq \sum_{1\leq i\leq m}\(\ve\|u_t\|_{\cC_q}^{q}+\ve^{-\b/(q-\b)}|u|^{\sig}_\sig\)
=m\ve\|u_t\|_{\cC_q}^{q}+mC_\ve|u|^{\sig}_\sig,
\end{split}
\end{equation*}
 where $\sig={q(q-1)}/({q-{\b}})$,  $\cC_q=\cC_{L^q(\W)}$, and $C_\ve=\ve^{-\b/(q-\b)}$. 
 One trivially verifies that  $\sig< {q_\gam}$.
Thus by Lemma \ref{l:2.0c} we deduce that
$$
C_\ve|u|^{\sig}_\sig\leq \de C_\ve |u|_{{q_\gam}}^{{q_\gam}}+\de^{-\sig/({q_\gam}-\sig)}C_\ve|\W|
$$
for any $\de>0$. Taking $\de =\ve /C_\ve$ in the above estimate, it gives
$$
C_\ve|u|^{\sig}_\sig\leq \ve |u|_{{q_\gam}}^{{q_\gam}}+\ve^{-q'}|\W|,\hs \mb{where }\,q'=\frac{\sig q+\b({q_\gam}-\sig)}{(q-\b)({q_\gam}-\sig)}.
$$
Therefore
$$
\int_\W|u|^{q-1}\,\max_{1\leq i\leq m}|u(t-r_i)|^{{\b}}dx\leq m\ve\|u_t\|_{\cC_q}^{q}+m\ve|u|_{{q_\gam}}^{{q_\gam}}+m\ve^{-q'}|\W|.
$$
It is trivial to check that
\be\label{eq'}
\lim_{q\ra\8}(q'/q)=1/(\gam-\b).\ee

We also infer from Lemma \ref{l:2.0c} that
$$
|u|_{q-2}^{q-2}\leq \ve |u|_{{q_\gam}}^{{q_\gam}}+\ve^{-(q-2)/(\gam+1)}|\W|,\hs |u|_{q-1}^{q-1}\leq \ve |u|_{{q_\gam}}^{{q_\gam}}+\ve^{-(q-1)/\gam}|\W|,
$$
and
$$
\int_\W |u|^{q-1}|h| dx\leq \ve|u|_{{q_\gam}}^{{q_\gam}}+\ve^{-(q-1)/\gam} |h|_{q_\gam/\gam}^{q_\gam/\gam}.$$
Combining all the above estimates together we obtain that
 \be\label{ode2}\begin{split}
 &\hs\frac{1}{q}\frac{d}{dt}|u|_{q}^{q}+(q-1)\int_\W|u|^{q-2}|\nab u|^2dx\\[1ex]
 &\leq -\(\Lam-\ve{(N+b_1'(m+1)+1)}\)|u|_{{q_\gam}}^{{q_\gam}}+{\ve m b_1'}\|u_t\|_{\cC_q}^{q}+C_\ve',
 \end{split}
 \ee
where
\be\label{e3.5}\begin{split}
C_\ve'=&\(N\ve^{-(q-2)/(\gam+1)}+b_1'(\ve^{-q'}+\ve^{-(q-1)/\gam})\)|\W|\\
&+\ve^{-(q-1)/\gam}\|h\|_{L^\8\(\R;L^{q_\gam/\gam}(\W)\)}^{q_\gam/\gam}\,.
\end{split}\ee

As ${q_\gam}>q$, by Lemma \ref{l:2.0c} one easily deduces  that
$$
|u|_{{q_\gam}}^{{q_\gam}}\geq |u|_{q}^{q}-|\W|.
$$
 We may assume $\Lam/2\geq \ve{(N+b_1'(m+1)+1)}$.  \eqref{ode2} then implies
\be\label{20b}\begin{split}
&\hs\frac{d}{dt}|u|_{q}^{q}+q(q-1)\int_\W|u|^{q-2}|\nab u|^2dx\\
&\leq -a_\ve |u|_{q}^{q}+{\ve q m b_1'}\|u_t\|_{\cC_q}^{q}+C_\ve'',
\end{split}  \ee
where
\be\label{e3.5b}a_\ve:=q(\Lam-\ve{(N+b_1'(m+1)+1)}),\hs C_\ve''=C_\ve'+a_\ve|\W|.\ee In particular,
\be\label{20c}\begin{split}
\frac{d}{dt}|u|_{q}^{q}\leq -a_\ve |u|_{q}^{q}+{\ve qm b_1'}\|u_t\|_{\cC_q}^{q}+C_\ve,
\end{split}  \ee

\vs
Let $E_\ve(t,s)=e^{-a_\ve(t-s)}$. Clearly for fixed $\ve>0$,
$$\lim_{t\ra +\8}E_\ve(t+s,s)=0$$
 uniformly w.r.t   $s\in\R$.
Multiplying  \eqref{20c} with $E_\ve(t,\tau)$ and integrating   in $t$ between $\tau$ and  $t$, it gives
\be\label{e:5.6}
y(t)\leq E_\ve(t,\tau)\|y_\tau\|_{C([-r,0])}+\int_\tau^t K_\ve(t,s)\|y_s\|_{C([-r,0])}ds+{C_\ve''}/{a_\ve}
\ee
for all $t\geq\tau\geq0$, where $y(t)=|u(t)|_{q}^{q}$, and $K_\ve(t,s)=\ve qmb_1' E_\ve(t,s)$.

\vs Note that $\vth:=\sup_{t\geq s}E_\ve(t,s)=1$. Since
$$a_\ve=q(\Lam-\ve{(N+b_1'(m+1)+1)})\geq q\Lam/2,$$
we have
\begin{equation*}\begin{split}
\int_0^tK_\ve(t,s)ds&=\ve qm b_1'\int_0^t E_\ve(t,s)ds\leq \ve qm b_1'\int_0^t e^{-q\Lam (t-s)/2}ds\\
&\leq \ve qm b_1'\int_0^\8 e^{-q\Lam t/2}dt\leq 2\ve m b_1'/\Lam.
\end{split}
\end{equation*}

Now we take
\be\label{e:e0}
\ve=\ve_0:=\frac{1}{2}\min\(\frac{\Lam}{2(N+b_1'(m+1)+1)}\,,\,\frac{\Lam}{4m b_1'}\).
\ee
Then  $a_{\ve_0}>q\Lam/2$ and
$$\ba{ll}\sup_{t\geq 0}\(\int_0^tK_{\ve_0}(t,s)ds\)\leq 1/4:=\k<1/2.\ea
$$
Thus by Lemma \ref{l:2.3} we deduce that $|u|_q$ is bounded on $[-r,T_\phi)$. Further using the same argument as in the proof of Theorem \ref{t:3.5} in Section 3.2 below with minor modifications, it can be shown that $|\nab u|$ is bounded on $[-r,T_\phi)$. It then follows that $T_\phi=\8$.

Clearly $
C_{\ve_0}''/a_{\ve_0}\leq 2C_{\ve_0}''/q\Lam:=\rho_q.
$
Thus by \eqref{e:5.6} we have
\be\label{e:3.8}
y(t)\leq E_{\ve_0}(t,\tau)\|y_\tau\|_{C([-r,0])}+\int_\tau^t K_{\ve_0}(t,s)\|y_s\|_{C([-r,0])}ds+\rho_q
\ee
for $t\geq\tau\geq0$. Now let us apply Lemma \ref{t:3.1} to the above inequality. The constants corresponding to those  in  Lemma \ref{t:3.1} and Remark \ref{r:2.5} read
$$
\mu=c=4/3,\hs \eta=7/2,\hs M=4/\sqrt 6,
$$
 and
$$
\lam=\lam_q:=\frac{\ln3-\ln2}{2(M_1+r a_{\ve_0})}\,a_{\ve_0},\hs\mb{where }\,M_1=\ln7-\ln2.
$$
By virtue of Lemma \ref{t:3.1} we conclude  that
\be\label{e:3.9}\ba{ll}
|u(t)|_{q}^{q}\leq Me^{-\lam_q t}\|\phi\|_{\cL^\8_q}^{q}+\frac{7}{2}\rho_q,\Hs \A\,t\geq 0.\ea
\ee
This completes the proof of the validity of \eqref{e:3.1}. $\Box$
\vs

 \bt\label{t:3.2}
 Suppose $h\in L^\8(\R\X\W).$ Then there exist $\rho_*,\lam_*>0$  such that
 \be\label{e:3.11a}\ba{lll}
|u(t;\phi)|_{\8}&\leq \left\{\ba{ll} \|\phi\|_{\cL^\8_\8}+\rho_*,\hs &\mb{if }\,r>0;\\[1ex]
e^{-\lam_*t}\|\phi\|_{\cL^\8_\8}+\rho_*,\hs &\mb{if }\,r=0
\ea\right.
\ea
\ee
for all $\phi\in \cC_{V_1}\Cap\cL^\8_\8$. 
   \et
{\bf Proof.} Let us first evaluate  $\limsup_{q\ra\8}\rho_q^{1/q}$. Note that
$q\Lam/2\leq a_{\ve_0}\leq q\Lam$. Therefore
$
\lim_{q\ra\8}a_{\ve_0}^{1/q}=1.
$
Hence
$$
\limsup_{q\ra\8}\rho_q^{1/q}=\limsup_{q\ra\8}\(C_{\ve_0}''/a_{\ve_0}\)^{1/q}=
\limsup_{q\ra\8}\(C_{\ve_0}''\)^{1/q}.
$$
On the other hand, by \eqref{e3.5}, \eqref{e3.5b}, Lemma \ref{l:2.0a} and \eqref{eq'}, we deduce that
\bes\begin{split}
\limsup_{q\ra\8}\(C_{\ve_0}''\)^{1/q}&\leq \limsup_{q\ra\8}2^{1/q}\(\(C_{\ve_0}'\)^{1/q}+(a_{\ve_0}|\W|)^{1/q} \)\\
&=\limsup_{q\ra\8}\(C_{\ve_0}'\)^{1/q}+1\leq \rho_*,
\end{split}
\ees
where
$$
\rho_*=\ve_0^{-{1}/{(\gam+1)}} + \ve_0^{-{1}/{(\gam-\b)}}+\ve_0^{-{1}/{\gam}}+\ve_0^{-{1}/{\gam}}\|h\|_{L^\8(\R\X\W)}^{1/\gam}+1
$$
Thus we conclude that  $\limsup_{q\ra\8}\rho_q^{1/q}\leq \rho_*$.

It is trivial to check  that if $r>0$ then $\limsup_{q\ra\8}\lam_q/q=0$; and if $r=0$, by the choice of $\ve_0$ we have
$$
\lam_q/q\geq \Lam(\ln3-\ln2)/(4M_1):=\lam_*
$$
for all $q$ sufficiently large. Therefore  using Lemma \ref{l:2.0a} once again,  we arrive  by \eqref{e:3.9} at the following estimates:
\bes\label{e:3.11}\ba{lll}
|u(t;\phi)|_{\8}&\leq \limsup_{q\ra\8}|u(t;\phi)|_{q}\\
&\leq \limsup_{q\ra\8}\(Me^{-\lam_q t}\|\phi\|_{\cL^\8_q}^{q}+\frac{7}{2}\rho_q\)^{1/q}\\[1ex]
&\leq \left\{\ba{ll} \|\phi\|_{\cL^\8_\8}+\rho_*,\hs &\mb{if }\,r>0;\\[1ex]
e^{-\lam_*t}\|\phi\|_{\cL^\8_\8}+\rho_*,\hs &\mb{if }\,r=0.
\ea\right.
\ea
\ees
This is precisely what we desired. $\Box$

\br\label{r:3.4}In  case   $r=0$ (i.e. for the equation without delay), we infer from Theorem \ref{t:3.2} that
$$
\limsup_{t\ra\8}|u(t;\phi)|_{\8}\leq \rho_*
$$
for all $\phi\in\cC_{V_1}\Cap\cL^\8_\8$. However, this remains an open problem in  case $r>0$, which may indicate some inherent differences between delay differential equations and those without delays. Fortunately, the following eventual  invariance property still remains true.
\er

\bp\label{p:3.4}Assume the hypotheses in Theorem \ref{t:3.2}, and let $\rho_*$ be the constant given therein. Then for any $\ve>0$, there is $t_0>0$ such that
\be\label{e:3.14}
|u(t;\phi)|_{\8}\leq \rho_*+\ve,\Hs \A\,t\geq t_0
\ee
for all $\phi\in\cC_{V_1}\Cap\cL^\8_\8$ with $\|\phi\|_{\cL^\8_\8}\leq \rho_*+\ve$.
\ep
{\bf Proof.}  Given $\ve>0$, to prove \eqref{e:3.14}, it suffices to check that for any $\rho<\rho_*+\ve:=\rho_*'$, the estimate holds true for each $\phi\in \cC_{V_1}\Cap\cL^\8_\8$ with $\|\phi\|_{\cL^\8_\8}\leq \rho$.

First,  it is easy to see that
\be\label{e:3.12}
\lim_{q\ra \8}\lam_q=r^{-1}(\ln3-\ln2):=2c_r.
\ee
Hence there is $q_0>0$ such that $\lam_q\geq c_r$ for all $q>q_0$. Take a $t_0>0$ such that $Me^{-c_r t_0}<7/2$. Then
\be\label{e:3.16}Me^{-\lam_q t}<7/2,\Hs \A\,t\geq t_0,\,\,q>q_0.\ee

Since $\limsup_{q\ra\8}\rho_q^{1/q}\leq \rho_*$,  there exists   $q_1>q_0$ such that
\be\label{e:3.17}
\rho_q\leq (\rho_*')^q,\Hs q>q_1.
\ee
Let $\|\phi\|_{\cL^\8_\8}\leq \rho$. We observe that
$$
\|\phi\|_{\cL^\8_q}^{q}\leq \|\phi\|_{\cL^\8_\8}^q |\W|\leq \(\rho|\W|^{1/q}\)^q.
$$
Because $\rho<\rho_*'$ and $|\W|^{1/q}\ra 1$ as $q\ra\8$, we can pick a $q_2>q_1$ such that
$
\rho|\W|^{1/q}<\rho_*'$ for all $q>q_2.$
It then follows that
\be\label{e:3.18}
\|\phi\|_{\cL^\8_q}^{q}\leq  \(\rho_*'\)^q,\Hs q>q_2.
\ee

Combing the above estimate together,  it yields
\be\label{e:3.19}\ba{ll}
Me^{-\lam_q t}\|\phi\|_{\cL^\8_q}^{q}\leq \frac{7}{2}\(\rho_*'\)^q,\Hs \A\, t\geq t_0,\,\,q>q_2.
\ea
\ee
Therefore by \eqref{e:3.17} and \eqref{e:3.19} we deduce that
$$\ba{ll}
|u(t;\phi)|_{\8}&\leq \limsup_{q\ra\8}\(Me^{-\lam_q t}\|\phi\|_{\cL^\8_q}^{q}+\frac{7}{2}\rho_q\)^{1/q}\\
&\leq \limsup_{q\ra\8}\(Me^{-\lam_q t}\|\phi\|_{\cL^\8_q}^{q}+\frac{7}{2}(\rho_*')^q\)^{1/q}\\[1ex]
&\leq\limsup_{q\ra\8}\(\frac{7}{2}(\rho_*')^q+\frac{7}{2}(\rho_*')^q\)^{1/q}=\rho_*'.
\ea
$$
The proof of the lemma is complete. $\Box$

\subsection{$H^1$ decay estimates}
In this subsection we give a decay estimate  in $V_1=H^1_0(\W)$.
   \bt\label{t:3.5}
Let $q_*<q<\8$. Suppose
$h\in L^\8(\R,L^{q_\gam/\gam}(\W))$.  Then there exist $M_1,\lam_1,\rho_1>0$ such that for all $\phi\in \cC_{V_1}\Cap \cL^\8_q$,
    \be\label{e3.21}\begin{split}
|\nab u(t;\phi)|^{2}&\leq \|\phi\|_{\cC_{V_1}}^2e^{-\mu_1 t}+M_1\|\phi\|_{\cL^\8_q}^{q}e^{-\lam_1 t}+\rho_1,\hs \A\,t\geq 0,\,\,\tau\in\R.
\end{split}
 \ee
   \et
{\bf Proof.} Let $u=u(t;\phi)$. Taking the inner product of  \eqref{eae} in $H$ with  $-\De u$, we obtain that
\be\label{e:3.23}\begin{split}
\frac{1}{2}\frac{d}{dt}|\nab u|^{2}+|\De u|^{2}&=-\(f(u)+G(t,u_t)+h(t),\,\De u\)\\[1ex]
 &\leq \frac{1}{2}|\De u|^{2}+ \frac{1}{2}|f(u)+G(t,u_t)+h(t)|^2\\[1ex]
 &\leq \frac{1}{2}|\De u|^{2}+ \frac{3}{2}\(|f(u)|^2+|G(t,u_t)|^2+|h(t)|^2\).
 \end{split}
 \ee
 By \eqref{e1.4}, \eqref{e1.7} and the H\"{o}lder's inequality, we deduce that
 \bes
 |f(u)|^2\leq a_1^2\int_\W(|u|^\a+1)^2dx\leq 2a_1^2\int_\W(|u|^{2\a}+1)dx\leq C_1|u|_q^{2\a}+C_2,
  \ees
  and
   \bes\begin{split}
  |G(t,u_t)|^2&\leq b_1^2\int_\W\(\,\max_{1\leq i\leq m}|u(t-r_i)|^{{\b}}+1\)^2dx\\[1ex]
  &\leq 2b_1^2\int_\W\(\,\max_{1\leq i\leq m}|u(t-r_i)|^{{2\b}}+1\)dx\\[1ex]
  &\leq2 b_1^2\int_\W\(\,\sum_{1\leq i\leq m}|u(t-r_i)|^{{2\b}}+1\)dx\leq C_3\|u_t\|_{\cC_q}^{2\b}+C_4.
 \end{split} \ees
Combing these estimates with \eqref{e:3.23} we find that
\be\label{e:3.24}
\frac{d}{dt}|\nab u|^{2}+|\De u|^{2}\leq 2C_1|u|_q^{2\a}+2 C_3\|u_t\|_{\cC_q}^{2\b}+C_5.
\ee
 Hence by Theorem \ref{t:3.1} one concludes that there exist $M',\lam'>0$ such that
\be\label{e:3.20}\begin{split}
\frac{d}{dt}|\nab u|^{2}&\leq -|\De u|^{2}+M'\|\phi\|_{\cL^\8_q}^qe^{-\lam' t}+C_6\\[1ex]
&\leq -\mu_1|\nab u|^{2}+M'\|\phi\|_{\cL^\8_q}^qe^{-\lam' t}+C_6,
 \end{split}
 \ee
 where $\mu_1$ is the first eigenvalue of $A=-\De$.

We may assume $\lam'\leq \mu_1/2$.
Then by the classical Gronwall lemma, there exists $\rho_1>0$ such   that
\be\label{e:3.25}\begin{split}
|\nab u|^{2}&\leq \|\phi\|_{\cC_{V_1}}^2e^{-\mu_1 t}+\frac{M'}{\mu_1-\lam'}\|\phi\|_{\cL^\8_q}^{q}\(e^{-\lam' t}-e^{-\mu_1 t}\)+\rho_1\\[1ex]
&\leq \|\phi\|_{\cC_{V_1}}^2e^{-\mu_1 t}+M_1\|\phi\|_{\cL^\8_q}^{q}e^{-\lam' t}+\rho_1,\Hs t\geq 0.
\end{split}
 \ee
 This verifies \eqref{e3.21}. $\Box$
 \br\label{r:3.6}
 Note that for $t\in [-r,0]$, since $u(t)=\phi(t)$ and $e^{-\mu_1 t}\geq 1$, we see that
  \eqref{e:3.25} readily  holds true. Therefore we actually have
 \be\label{e:e1}\begin{split}
|\nab u|^{2}&\leq \|\phi\|_{\cC_{V_1}}^2e^{-\mu_1 t}+M_1\|\phi\|_{\cL^\8_q}^{q}e^{-\lam' t}+\rho_1:=e_1(t),\Hs t\geq -r.
\end{split}
 \ee
 \er
 \br\label{r:3.7}
Integrating  \eqref{e:3.20}  between $t$ and $t+1$,  by \eqref{e:e1} one obtains that
\be\label{e:3.26}\begin{split}
\int_t^{t+1}|\De u|^2ds&\leq \(|\nab u(t)|^{2}+|\nab u(t+1)|^{2}\)+ \frac{M'}{\lam'}\|\phi\|_{\cL^\8_q}^qe^{-\lam' t} +C_6  \\[1ex]
&\leq 2e_1(t)+\frac{M'}{\lam'}\|\phi\|_{\cL^\8_q}^qe^{-\lam' t} +C_6\\[1ex]
&\leq C_7\(\|\phi\|_{\cC_{V_1}}^2+2\|\phi\|_{\cL^\8_q}^{q}\)e^{-\lam' t}+C_8:=e_2(t).
\end{split}
\ee

Similarly, if we integrate \eqref{20b} with $\ve=\ve_0$ ($\ve_0$ is given by \eqref{e:e0}) between $t$ and $t+1$,   by \eqref{e:3.1} we get
\be\label{e:3.28}\begin{split}
\int_t^{t+1}\int_\W|u|^{q-2}|\nab u|^2dx\,ds&\leq \(|u(t)|_{q}^{q}+|u(t+1)|_{q}^{q}\)\\
&\hs +\ve_0qmb_1'\int_t^{t+1}\|u_s\|_{\cC_q}^qds+C_{\ve_0}''\\[.5ex]
&\leq C_9\|\phi\|_{\cL^\8_q}^{q}e^{-\lam_q t}+C_{10}:=e_3(t).
\end{split}
\ee
 For  $T>0$, integrating  \eqref{20b} with $\ve=\ve_0$ between $0$ and $T$,   it also yields
\be\label{e:3.28b}\begin{split}
\int_0^T\int_\W|u|^{q-2}|\nab u|^2dx\,ds\leq C_T\(\|\phi\|_{\cL^\8_q}^{q}+1\)
\end{split}
\ee
for some $C_T>0$.
 \er

 \section{$H^2$ estimates}
 This section is devoted to the $H^2$ decay estimates and global estimates of solutions of \eqref{e:ivp2}.
 As in Section3, we may assume the initial time $\tau=0$. This is just for the sake of simplicity in writing, and all the estimates remain true for \eqref{e:ivp2} uniformly with respect to $\tau\in\R$.
\subsection{The case $\phi\in \cL^\8_\8$}
We first give an $H^2$ estimate in case the initial value $\phi\in  \cL^\8_\8$.

   \bp\label{p:3.7}
 Suppose
$h\in L^\8(\R\X\W)\Cap L^\8(\R;H^1(\W))$.  Then for any $R>0$, there  exist $C=C(R)>0$ such that
 \be\label{e3.21}\begin{split}
|\De u(t;\phi)|\leq C,\Hs  t\geq r+1
\end{split}\ee
 for all $\phi\in \cC_{V_1}\Cap \cL^\8_\8$ with
 \be\label{e:3.32}\|\phi\|_{\cC_{V_1}},\,\|\phi\|_{\cL_\8^\8}\leq R.
     \ee
   \ep
{\bf Proof.}  Let $u=u(t;\phi)$. Multiplying \eqref{eae} with $-\De u'$ (where \, $'=\frac{\pa}{\pa t}$) and integrating over $\W$,  we get
\be\label{e:3.33}\begin{split}
|\nab u'|^2+\frac{1}{2}\frac{d}{dt}|\De u|^2= -\int_\W\(f(u)+G(t,u_t)+h\)\De  u'dx= I_1+I_2+I_3,
\end{split}
 \ee
 where $I_3=\int_\W \nab h\.\nab  u'dx$, and
 $$
I_1= \int_\W f'(x,u)\nab u\.\nab  u'dx,
$$
$$
I_2=\int_\W\(\sum_{i=1}^m \pa_i g(u(t-r_1),\cdots,u(t-r_m))\nab u(t-r_i)\)\.\nab  u'dx.$$
 Here and below $\pa_i g(v)=\frac{\pa }{\pa v_i}g(v)$ ($v\in\R^m$).

\vs Let us first evaluate $I_1$. Fix a $q_*< q<\8$. Using  the structure conditions (F1) and (G1),  the H\"{o}lder's inequality and the Young's inequality we deduce that
\be\label{e:3.29}\begin{split}
I_1&\leq \int_\W|f'(x,u)||\nab u|\,|\nab u'|dx\\
&\leq a_0 \int_\W(|u|^{\a-1}+1)|\nab u|\,|\nab u'|dx+ \\
&\leq \frac{1}{4}|\nab u'|^2+a_0^2\int_\W\(|u|^{\a-1}|\nab u|+1\)^2dx\\
&\leq \frac{1}{4}|\nab u'|^2+2a_0^2\int_\W|u|^{2(\a-1)}|\nab u|^2dx+2a_0^2|\W|\\
&\leq \frac{1}{4}|\nab u'|^2+2a_0^2\int_\W(|u|^{q-2}+1)|\nab u|^2dx+2a_0^2|\W|.
\end{split}
 \ee
 Here we have used the simple fact $|u|^{2(\a-1)}<|u|^{q-2}+1$.
Now assume that $\phi\in \cC_{V_1}\Cap \cL^\8_\8$ and satisfies \eqref{e:3.32}.
  Then by Theorem \ref{t:3.2} there exists $\rho_*>0$ (independent of $R$) such that
\be\label{e:3.30}
|u(t)|_{\8}\leq  R+\rho_*,\Hs t\geq0.
\ee
Hence
\bes\label{e:3.31}\begin{split}
I_1 &\leq \frac{1}{4}|\nab u'|^2+C_1(R)|\nab u|^2+C_{12}\\[.5ex]
&\leq \frac{1}{4}|\nab u'|^2+ C_1(R)e_1(t)+C_{12},\Hs t\geq 0,
\end{split}
 \ees
where $e_1(t)$ is the function defined in \eqref{e:e1}.

\vs
 Note that
$$
I_2\leq \int_\W|\nab g(u(t-r_1),\cdots,u(t-r_m))|\(\sum_{i=1}^m |\nab u(t-r_i)|\)\,|\nab u'|dx.
$$
By \eqref{e1.6} we deduce that
\bes\begin{split}
|\nab g(u(t-r_1),\cdots,u(t-r_m))|&\leq b_0'\(\max_{1\leq i\leq m}|u(t-r_i)|^{\b-1}+1\)\\
&\leq b_0' \(\sum_{1\leq i\leq m}|u(t-r_i)|^{\b-1}+1\).
\end{split}\ees
Thus a similar argument as above applies to show  that  
 \be\begin{split}
I_2&\leq b_0' \int_\W\(\sum_{1\leq i\leq m}|u(t-r_i)|^{\b-1}+1\)\(\sum_{i=1}^m |\nab u(t-r_i)|\)\,|\nab u'|dx\\
&\leq \frac{1}{4}|\nab u'|^2+C_2(R)\sum_{i=1}^m|\nab u(t-r_i)|^2\\[.5ex]
&\leq(\mb{by Remark }\ref{r:3.6})\leq \frac{1}{4}|\nab u'|^2+C_2(R)\sum_{i=1}^m e_1(t-r_i)\\
&\leq \frac{1}{4}|\nab u'|^2+C_3(R)e_1(t-r),\Hs t\geq 0.
\end{split}
 \ee

We also have
\be\label{e:3.38}
 I_3\leq \frac{1}{4}|\nab u'|^2+|\nab h|^2\leq \frac{1}{4}|\nab u'|^2+\|h\|_{L^\8(\R;\,H^1(\W))}^2,\Hs t\geq 0.
\ee
 Therefore by \eqref{e:3.33} one concludes that
 \be\label{e:3.34}
 \begin{split}
 \frac{1}{4}|\nab u'|^2+ \frac{1}{2}\frac{d}{dt}|\De u|^2&\leq  C_1(R)e_1(t)+C_3(R)e_1(t-r)+C_{13}\\
 &\leq C_4(R)e_1(t-r)+C_{13}\\
 &:=a(t),\Hs t\geq 0.
 \end{split}
 \ee

Since $e_i(t)$ ($i=1,2,3$) are nonincreasing,  we have
$$
\int_t^{t+1}a(s)ds\leq a(t)\leq a(t_0),\Hs \A\,t\geq t_0\geq r,
$$
$$
\int_t^{t+1}|\De u|^2dt\leq (\mb{by }\eqref{e:3.26})\leq  e_2(t_0), \Hs \A\,t\geq t_0\geq r.
$$
Applying the Uniform Gronwall Lemma (see Temam \cite[pp. 89, Lemma 1.1]{Tem}) to \eqref{e:3.34} yields
$$
|\De u(t+1)|^2\leq 2a(t_0)+e_2(t_0):=C_5(R),\Hs \A\, t\geq t_0\geq r.
$$
Taking $t_0=r$ one immediately obtains \eqref{e3.21}. $\Box$

\bt\label{t:3.8}  Suppose
$h\in L^\8(\R\X\W)\Cap L^\8(\R;H^1(\W))$.  Then for any $R>0$, there  exist $C=C(R)>0$ such that
 \be\label{e3.21b}\begin{split}
|\De u(t;\phi)|\leq C,\Hs t\geq 0
\end{split}\ee
 for all  $\phi\in \cC_{V_2}\Cap \cL^\8_\8$ with
 $\|\phi\|_{\cC_{V_2}},\,\|\phi\|_{\cL_\8^\8}\leq R.$
   \et
{\bf Proof.} By \eqref{e:3.34}  one easily deduces that there is a constant  $C=C(R)>0$ such that
$|\De u(t;\phi)|\leq C$ for $t\in[0,r+1]$. \eqref{e3.21b} then directly follows from this local estimate and Proposition \ref{p:3.7}. $\Box$

\subsection{The case of separated delays}
In this part we consider a slightly particular case where the delays $r_i$ ($1\leq i\leq m$) are separated, namely, $g$ takes the form
\be\label{e:3.35}
g(v)=g_1(v_1)+g_2(v_2)+\cdots+g_m(v_m),\Hs v\in\R^m.
\ee
In such a case we establish  some $H^2$ decay estimates by assuming $\phi\in \cL^\8_q$ for some $q> q_*$ rather than  $\phi\in \cL^\8_\8$.

We begin with \eqref{e:3.33}. By \eqref{e:3.29} and \eqref{e:e1} we have
\bes\label{e:3.39}\begin{split}
I_1
&\leq \frac{1}{4}|\nab u'|^2+ 2a_0^2\int_\W|u|^{q-2}|\nab u|^2dx+2a_0^2e_1(t)+2a_0^2|\W|.
\end{split}
 \ees
Using the H\"{o}lder's inequality, the Chauchy-Schwartz inequality, the structure condition (G1) and \eqref{e:e1}, we deduce that
 \bes\label{e:3.40}\begin{split}
I_2&=\sum_{i=1}^m \int_\W g_i'(u(t-r_i))\nab u(t-r_i)\.\nab u'\,dx\\
&\leq \frac{1}{4}|\nab u'|^2+m\sum_{i=1}^m\int_\W|g_i'(u(t-r_i))|^2|\nab u(t-r_i)|^2dx\\[.5ex]
&\leq \frac{1}{4}|\nab u'|^2+C_{11}\sum_{i=1}^m\int_\W \(|u(t-r_i)|^{2(\b-1)}+1\)|\nab u(t-r_i)|^2dx\\
&\leq \frac{1}{4}|\nab u'|^2+C_{11}\sum_{i=1}^m\int_\W \(|u(t-r_i)|^{q-2}+2\)|\nab u(t-r_i)|^2dx\\
&\leq \frac{1}{4}|\nab u'|^2+C_{11}\sum_{i=1}^m\int_\W |u(t-r_i)|^{q-2}|\nab u(t-r_i)|^2dx+C_{12}\sum_{i=1}^m e_1(t-r_i)\\
&\leq \frac{1}{4}|\nab u'|^2+C_{11}\sum_{i=1}^m\int_\W |u(t-r_i)|^{q-2}|\nab u(t-r_i)|^2dx+C_{13}\,e_1(t-r),\\
\end{split}
 \ees
where $C_{11}=2m b_0^2$. Combining \eqref{e:3.33}, \eqref{e:3.38} and the above two estimates  it yields
\be\label{e:3.41}
 \begin{split}
 \frac{1}{4}|\nab u'|^2+ \frac{1}{2}\frac{d}{dt}|\De u|^2&\leq  2a_0^2\int_\W|u|^{q-2}|\nab u|^2dx\\
 &\hs+C_{11}\sum_{i=1}^m\int_\W |u(t-r_i)|^{q-2}|\nab u(t-r_i)|^2dx\\
 &\hs +C_{14} \,e_1(t-r)+C_{15}:=I(t),\Hs t\geq 0.
 \end{split}
 \ee

By Remark \ref{r:3.7} we find that
\bes\label{e:3.45}
 \begin{split}
\int_t^{t+1}I(t)dt&\leq  2a_0^2\,e_3(t_0)+C_{11}\sum_{i=1}^m\,e_3(t_0-r_i)+C_{14} \,e_1(t_0-r)+C_{15}\\
&\leq C_{16}\,e_3(t_0-r)+C_{14} \,e_1(t_0-r)+C_{15},\Hs \A\,t\geq t_0\geq r,
\end{split}
 \ees
$$
\int_t^{t+1}|\De u|^2dt\leq e_2(t_0), \Hs \A\,t\geq t_0\geq r.
$$
Thanks to the Uniform Gronwall Lemma, 
it follows from  \eqref{e:3.41} and the definition of $e_1(t),e_2(t)$ and $e_3(t)$ (see \eqref{e:e1} and Remark \ref{r:3.7}) that
\bes
 \begin{split}
|\De u(t+1)|^2&\leq 2e_2(t_0)+2\(C_{16}\,e_3(t_0-r)+C_{14} \,e_1(t_0-r)+C_{15}\)\\
&\leq M_2'\(\|\phi\|_{\cC_{V_1}}^2+\|\phi\|_{\cL_q^\8}^q\)e^{-\lam_2(t_0-r)}+\rho_2,\Hs t\geq t_0
\end{split}\ees
for some constants $M_2',\lam_2,\rho_2>0$. In particular,
\be\label{e:3.43}
 \begin{split}
|\De u(t_0+1)|^2\leq M_2\(\|\phi\|_{\cC_{V_1}}^2+\|\phi\|_{\cL_q^\8}^q\)e^{-\lam_2(t_0+1)}+\rho_2.
\end{split}\ee
where $M_2=M_2'e^{\lam_2(r+1)}$. Rewriting $t_0+1=t$ in \eqref{e:3.43}, we get

   \bp\label{p:3.10}
Let $q_*<q<\8$. Suppose
$h\in L^\8(\R\X\W)\Cap L^\8(\R;H^1(\W))$.  Then there exist $M_2,\lam_2>0$ and $\rho_2>0$ such that for all  $\phi\in \cC_{V_1}\Cap \cL^\8_q$,
\be\label{e:3.44}
 \begin{split}
|\De u(t;\phi)|^2\leq M_2\(\|\phi\|_{\cC_{V_1}}^2+\|\phi\|_{\cL_q^\8}^q\)e^{-\lam_2t}+\rho_2,\Hs t\geq r+1.
\end{split}\ee
   \ep

   In addition to the hypotheses in Proposition \ref{p:3.10},
   if we assume
   \be\label{e:3.47}\int_{-r}^0\int_\W|\phi|^{q-2}|\nab \phi|^2dx\,dt<\8,\ee
    then for $T=r+1$, by \eqref{e:3.28b}  and the definition of $e_1(t)$ it is easy to deduce that
   $$
\int_0^TI(t)dt\leq C_{17} \(\|\phi\|_{\cC_{V_1}}^2+\|\phi\|_{\cL_q^\8}^q+\int_\W|\phi|^{q-2}|\nab \phi|^2dx\) +C_{18}. $$
Integrating \eqref{e:3.41} between $0$ and $T$ one finds  that
$$
|\De u(t)|^2\leq \|\phi\|^2_{\cC_{V_2}}+2C_{17} \(\|\phi\|_{\cC_{V_1}}^2+\|\phi\|_{\cL_q^\8}^q+\int_\W|\phi|^{q-2}|\nab \phi|^2dx\) +2C_{18}
$$
for $ t\in[0,T]$. Combining this with Proposition \ref{p:3.10} it yields   the following result.

 \bt\label{t:3.11}Assume $g$ takes the form in \eqref{e:3.35}.
Let $q_*<q<\8$. Suppose
$h\in L^\8(\R\X\W)\Cap L^\8(\R;H^1(\W))$.  Then there exist positive constants $\~M_2$ and $\~\rho_2$ such that
\bes\label{e:3.46}
 \begin{split}
|\De u(t;\phi)|^2\leq \~M_2\(\|\phi\|_{\cC_{V_2}}^2+\|\phi\|_{\cL_q^\8}^q+\int_\W|\phi|^{q-2}|\nab \phi|^2dx\)e^{-\lam_2t}+\~\rho_2,\Hs t\geq 0
\end{split}\ees
 for all  $\phi\in \cC_{V_2}\Cap \cL^\8_q$ satisfying \eqref{e:3.47}.
   \et

\section{Existence, Uniqueness and Regularity of Solutions}
 Using the estimates established in the previous sections, it can be shown
  by very standard argument via Galerkin approximation methods as stated in the beginning of Section 3 that the following existence and uniqueness result hold true for the initial value problem \eqref{e:ivp2}.

 \bt\label{t:5.1} Let $q_*<q\leq \8$. Suppose
$h\in L^\8(\R;L^{q_\gam/\gam}(\W))$.  Then for each $\phi\in \cC_{V_1}\Cap \cL^\8_q$, problem \eqref{e:ivp2} has a unique global weak solution $u=u(t;\tau,\phi)$ (in the distribution sense) with
\be\label{e:5.1}
u\in C\([-r,\8);V_1\)\Cap L^\8\([-r,\8);V_1\)\Cap L^\8\([-r,\8);L^q(\W)\).
\ee
Furthermore, for any $0<T<\8$,
\be\label{e:5.2}
u\in L^2\([0,T];V_2\).
\ee
   \et
   \br Since $u\in L^\8\([-r,\8);L^q(\W)\)$ and $q>q_*$, by \eqref{e1.4} and \eqref{e1.5} one trivially verifies
   that
   \be\label{e:5.3}w(t):=f(u)+G^\tau(t,u_t)+h^\tau(t)\ee belongs to $L^2([0,T];H)$ for any $T>0$. The  relation $u\in C\([-r,\8);V_1\)$ in \eqref{e:5.1} then follows from Theorem 3.3 in \cite[Chapt. II]{Tem} on abstract linear equations.
   \er

Now we pay some more attention to regularity of solutions when the initial value $\phi$ has higher regularity.
Let $w$ be the function in \eqref{e:5.3}. For simplicity, as before we set $\tau=0$, hence $h^\tau(t)=h(t)$, and
$$
G^\tau(t,u_t)=G(t,u_t)=g(u(t-r_1),\cdots,u(t-r_m)).
$$
Therefore
\be\label{e:w'}
w'=f'(u)u'+\sum_{1\leq i\leq m}\pa_i g(u(t-r_1),\cdots,u(t-r_m))u'(t-r_i)r_i'+h'.
\ee
As above, one easily verify that all the functions $f'(u)$ and $\pa_i g(u(t-r_1),\cdots,u(t-r_m))$ ($1\leq i\leq m$) belong to $L^2([0,T];H)$ for any $T>0$.

 Multiplying \eqref{eae} with $u'$ and integrating over $\W$,  we get
\be\label{e:5.4}\begin{split}
|u'|^2+\frac{1}{2}\frac{d}{dt}|\nab u|^2&= \(f(u)+G(t,u_t)+h,\, u'\)\\
&\leq \frac{1}{2}|u'|^2+\frac{1}{2}|f(u)+G(t,u_t)+h|^2.
\end{split}
 \ee
 Using the above inequality and the $H^1$ estimates it is not difficulty to see that $u'\in L^2((0,T);H)$ for any $T>0$. Further if  we assume that $\phi'\in L^2((0,T);H)$ then since $u\in C([-r,T];H)$, it can be shown that $u'\in L^2((-r,T);H)$.

Now assume  $r_i\in C^1(\R)$. Then by \eqref{e:w'} we deduce that $w'\in L^2((0,T);H)$.
 Thanks to Theorem 3.2 in \cite[Chapt. II]{Tem} on regularity of abstract linear equations, and the $H^2$ estimates given in Section 4, we obtain the following theorems.
 \bt\label{t:5.2}Assume $r_i\in C^1(\R)$ $(1\leq i\leq m)$, and that
$$
 h\in L^\8(\R\X\W)\Cap L^\8(\R;H^1(\W)),$$
 $$
  h'\in L^2((a,b);H),\Hs \A\,-\8<a<b<\8.$$
   Then  for any   $\phi\in \cC_{V_2}\Cap \cL^\8_\8$ with $\phi'\in L^2((-r,0);H)$, the solution $u$ of \eqref{e:ivp2} given by Theorem \ref{t:5.1} satisfies
   $$
   u\in C([-r,\8);V_2),
   $$
   $$ u'\in L^2((0,T);V_1)\Cap C([0,T];H),\Hs \A\,0<T<\8.
   $$
   \et

 \bt\label{t:5.3}Assume $g$ takes the form in \eqref{e:3.35}, and that $r_i\in C^1(\R)$ $(1\leq i\leq m)$. Suppose $$
 h\in L^\8(\R\X\W)\Cap L^\8(\R;H^1(\W)),$$
 $$
  h'\in L^2((a,b);H),\Hs \A\,-\8<a<b<\8.$$
   Then  for any   $\phi\in \cC_{V_2}\Cap \cL^\8_q$ $(q_*<q<\8)$ satisfying \eqref{e:3.47} and $\phi'\in L^2((-r,0);H)$, the solution $u$ of \eqref{e:ivp2} given by Theorem \ref{t:5.1} satisfies
   $$
   u\in C([-r,\8);V_2),
   $$
   $$ u'\in L^2((0,T);V_1)\Cap C([0,T];H),\Hs \A\,0<T<\8.
   $$
   \et

\Vs
\end{document}